\documentclass[12pt]{article}

\usepackage{amssymb}

\def\M{{\cal M}}
\def\oM{{\overline{\cal M}}}
\def\Q{{\mathbb Q}}
\def\A{{\cal A}}
\def\qed{{\hfill $\diamondsuit$}}
\def\cL{{\cal L}}
\def\CP{{{\mathbb C}{\rm P}}}
\def\Aut{{\rm Aut}}

\usepackage{theorem}

\newtheorem{theorem}{Theorem}
\newtheorem{proposition}{Proposition}[section]
\newtheorem{corollary}[proposition]{Corollary}

{\theorembodyfont{\rmfamily}
\newtheorem{definition}[proposition]{Definition}
\newtheorem{example}[proposition]{Example}
\newtheorem{remark}[proposition]{Remark}
\newtheorem{notation}[proposition]{Notation}
}

\include{epsf}

\title{Enumeration of ramified coverings of the sphere and 
2-dimensional gravity}

\author{Dimitri Zvonkine\thanks{
Institut math{\'e}matique de Jussieu,
Universit{\'e} Paris~VI, 175, rue du Chevaleret,
75013 Paris, France. E-mail: zvonkine@math.jussieu.fr\,.\/\/
The author was partially suported by
EAGER - European Algebraic Geometry Research Training Network, 
contract No. HPRN-CT-2000-00099 (BBW) and by the Russian Foundation
of Basic Research grant 02-01-22004.}}

\date{\today}

\begin{document}

\maketitle

\begin{abstract}
Let $\A$ be the algebra generated by the power series
$\sum n^{n-1} q^n/n!$ and $\sum n^n q^n /n!\,$. We prove
that many natural generating functions lie in this algebra:
those appearing in graph enumeration problems, in the 
intersection theory of moduli spaces $\oM_{g,n}$ and in the enumeration
of ramified coverings of the sphere. 

We argue that ramified coverings of the sphere with a large number
of sheets provide a model of 2-dimensional gravity.
Our results allow us to compute the asymptotic of the number
of coverings as the number of sheets goes to infinity.
The leading terms of such asymptotics are the values of certain
observables in 2-dimensional gravity. We prove that they coincide
with the values provided by other models. In particular, we 
recover a solution of the Painlev{\'e}~I equation and the string
solution of the KdV hierarchy.
\end{abstract}

\tableofcontents

\section{Introduction}

Denote by $\A$ the subalgebra of the algebra of
power series in one variable, generated by the series
$$
\sum_{n \geq 1} \frac{n^{n-1}}{n!} \, q^n \quad \mbox{and} \quad
\sum_{n \geq 1} \frac{n^n}{n!} \, q^n.
$$
We wish to show that this algebra plays an important role
in the intersection theory of moduli spaces $\oM_{g,n}$
of stable curves and in the problem of enumeration of
ramified coverings of the sphere.

\paragraph{SECTION~\ref{Sec:A}} 
contains a more explicit description of the algebra $\A$. 
In this section we also prove some relations between $\A$
and the combinatorics of Cayley trees ($=$ trees with
numbered vertices).

\paragraph{SECTION~\ref{Sec:coverings}}
is devoted to the problem of enumerating the ramified
coverings of the sphere with specified ramification types.

Consider a holomorphic map $f : C \rightarrow \CP^1$
of degree $n$ from a smooth complex curve $C$ to the Riemann sphere. 
Such maps will be called {\em ramified coverings} with $n$
sheets.

A {\em ramification point} of $f$ is a point of the target Riemann
sphere that has less than $n$ distinct preimages. 

For each ramification point $y$ of a ramified covering, we are
going to single out several {\em simple} preimages of $y$.

\begin{definition} \label{Def:markedcov}
A {\em marked ramified covering} is a ramified covering
with a choice, for every ramification point $y$, of a subset
of the set of simple preimages of $y$.
\end{definition}

Consider a partition $\mu = 1^{a_1} 2^{a_2} \dots$ of an integer 
$m \leq n$. (Here we use multiplicative notation for
partitions: the partition $\mu$ contains $a_1$ parts equal to
$1$, $a_2$ parts equal to $2$, and so on, $\sum i a_i = m$.)
Suppose that a point $y \in \CP^1$ has $a_1$ 
marked simple preimages, $a_2$ double preimages, and so on. 
(Consequently, $y$ also has $n-m$ unmarked simple preimages.)
We then say that $y$ is a ramification point of $f$ of {\em multiplicity}
$r=a_2+2a_3 + 3a_4+\dots$ and of {\em ramification type}
$\mu = 1^{a_1}2^{a_2} \dots$. Sometimes the number $r$ will
also be called the {\em degeneracy} of the partition $\mu$.

\begin{definition}
\label{Def:Hurwitz}
A {\em Hurwitz number} $h_{n; \mu_1, \dots, \mu_k}$ 
is the number of connected $n$-sheeted marked ramified coverings of $\CP^1$
with $k$ ramification points, whose ramification types
are $\mu_1, \dots, \mu_k$. Every such covering is counted with weight
$1/|\Aut|$, where $|\Aut|$ is the number of automorphisms of the
covering.
\end{definition}

Note that the genus $g$ of the covering surface can be reconstituted
from the data $(n; \mu_1, \dots, \mu_k)$ using the Riemann-Hurwitz
formula: if the degeneracy of $\mu_i$ equals $r_i$, then
$$
2-2g = 2n-r_1- \dots - r_k.
$$

\bigskip

Fix $k$ nonempty partitions $\mu_1, \dots, \mu_k$ with
degeneracies $r_1, \dots, r_k$. Let $r$
be the sum $r = r_1 + \dots + r_k$. 

\begin{notation}
\label{Not:Hurwitz}
Denote by $h_{g,n;\mu_1, \dots, \mu_k}$ the number of $n$-sheeted
marked ramified coverings of $\CP^1$ by a genus $g$ surface,
with $k$ ramification points of types $\mu_1, \dots, \mu_k$
and, in addition, $c(n) = 2n+2g-2-r$ simple ($=$ of multiplicity~1)
ramification points. Each covering is counted with weight
$1/|\Aut|$.
\end{notation}

\begin{theorem} \label{Thm:Hurwitz}
Fix any $g \geq 0$, $k \geq 0$. If $g=1$, we suppose that
$k \geq 1$. Then for any partitions
$\mu_1, \dots, \mu_k$, the series
$$
H_{g;\mu_1, \dots, \mu_k} (q) =
\sum_{n \geq 1} \frac{h_{g,n;\mu_1, \dots, \mu_k}}{c(n)!} 
\, q^n
$$
lies in the algebra $\A$.
\end{theorem}

The only known proof of this theorem involves a surprising
detour by the intersection theory on moduli spaces of stable
curves. Using the Ekedahl-Lando-Shapiro-Vainshtein (or ``ELSV'') 
formula, one can express the Hurwitz numbers for $k=1$ as
integrals of some cohomology classes over these moduli spaces.
In~\cite{GoJaVa} the ELSV formula is used to express the
generating functions for Hurwitz numbers in the case $k=1$
as rational functions of the series $Y(q)$. This essentially
proves the theorem for $k=1$. After that, we proceed by induction
on~$k$.

Among other things, the theorem allows one to find the
asymptotic of the coefficients of $H_{g;\mu_1, \dots, \mu_k}$ as 
$n \rightarrow \infty$, knowing only the several first
coefficients.

\paragraph{SECTION~\ref{Sec:gravity}} describes a
model of $2$-dimensional gravity obtained by counting 
ramified coverings of the sphere.
We show that one can extract from the symptotic of Hurwitz
numbers a solution of the Painlev{\'e}~I equation and
the string solution of the Korteweg~-~de~Vries equation.
The same solutions are obtained in other models
of $2$-dimensional gravity (by counting quadrangulations
or using intagrals over moduli spaces of curves).

We also compare the enumerative problems concerning ramified
coverings of the sphere and those of the torus. While the former
are related to the intersection theory on $\oM_{g,n}$
and give rise to the algebra $\A$, the latter are related
to volumes of spaces of abelian differentials on Riemann
surfaces and give rise to the algebra of quasi-modular forms.

\paragraph{Acknowledgments}

The author is grateful to J.-M.~Bismut, F.~Labourie,
M.~Kontsevich, S.~Natanzon, Ch.~Okonek,
A.~Okounkov, D.~Panov, J.-Y. Welschinger,
D.~Zagier, A.~Zorich and A.~Zvonkin for useful
discussions and remarks. A special thank to Sergei Lando,
with whom we proved together some of the results of the 
last section, and to M.~Kazarian for sharing his own work 
on the same subject. I would also like to thank for their interest 
the participants of the mathematical physics seminar at the ETH Z{\"u}rich
and of the mathematical seminar at the ENS Lyon, as
well as the participants of the Luminy conference on
billiards and Teichm{\"u}ller spaces.

This work was partially supported by
EAGER - European Algebraic Geometry Research Training Network, 
contract No. HPRN-CT-2000-00099 (BBW) and by the RFBR grant
02-01-22004.

\paragraph{Notation.} Here we summarize some
notation that we use consistently throughout the
paper.

\bigskip

\noindent
\begin{tabular}{lp{29em}}
$n$ & The number
of sheets of a covering.
The power of the variable $q$ in generating series. 
The number of marked points on a Riemann surface is
sometimes $n$ and sometimes $n-r$.\\
$q$ & The variable in generating series (to a sequence
$s_n$ we usually assign the series $\sum s_n q^n/n!$).\\
$g$ & The genus of a Riemann surface.\\
$\mu$ & A partition.\\
$p$ & The number of parts of a partition $\mu$.\\
$a_i$ & The number of parts of a partition $\mu$ that are
equal to $i$.\\
$b_i$ & The parts of a partition $\mu$
are denoted by $b_1, \dots, b_p$.\\
$r$ & The degeneracy of a partition $\mu$ defined
by $r= \sum (b_i-1)$. The multiplicity
of a ramification point.\\
$k$ & The number of partitions. If $k > 1$, the
partitions are denoted by $\mu_1, \dots, \mu_k$,
their degeneracies by $r_1, \dots, r_k$, while
$r$ is the total degeneracy $r = \sum r_i$.\\
$c(n)$ & The number of simple ramification points
in a ramified covering.\\
$\psi_i$ & The first Chern class $c_1(\cL_i)$
of the line bundle $\cL_i$ over $\oM_{g,n}$.\\
$d_i$ & The power of the class $\psi_i$ in the
intersection numbers we consider. 
\end{tabular}

\section{The algebra $\A$ of power series}
\label{Sec:A}

The algebra of power series
$$
\A = \Q \left[ \sum_{n \geq 1} \frac{n^{n-1}}{n!} \, q^n,
\sum_{n \geq 1} \frac{n^n}{n!} \, q^n \right]
$$
plays a central role in this paper.
Here we give an explicit description of $\A$ and show its
relation with the combinatorics of Cayley trees. Many of
the results below are known, but have probably never
been put together. As far as we know, the algebra $\A$
itself was first discovered by D.~Zagier several years ago
(unpublished), and then independently 
introduced in our paper~\cite{Zvonkine},
where most of the results of Section~\ref{Ssec:compA} are
given. Various series from $\A$ also appear in~\cite{GoJaVa}.

\subsection{How to make computations in $\A$}
\label{Ssec:compA}

Denote by $Y$ and $Z$ the generators of $\A$
$$
Y = \sum_{n \geq 1} \frac{n^{n-1}}{n!} \, q^n, \qquad
Z = \sum_{n \geq 1} \frac{n^n}{n!} \, q^n.
$$
Denote by $D$ the differential operator 
$D = q \, \frac{\partial}{\partial q}$. Thus $Z = DY$.

Note that both $Y$ and $Z$ have a radius of convergence
of $1/e$. Therefore the same is true of all series
in $\A$. The function $Y(q)$, more precisely, $-Y(-q)$,
was considered by J.~H.~Lambert~\cite{Lambert} in
1758\footnote{We thank N. A'Campo for this reference.}. The relations
that follow can be deduced from the Lagrange inversion
theorem applied to the equation $Y(q) = q e^{Y(q)}$ 
or from the Abel identities (see~\cite{GouJacbook}, Section~1.2).

\begin{proposition}\label{Prop:Y=qe^Y}
We have
$$
Y = q e^Y.
$$
\end{proposition}

\paragraph{Proof.} $Y$ is the exponential generating series
for rooted Cayley trees (Definition~\ref{Def:Cayley}). 
Therefore $e^Y$ is the exponential generating series for 
forests of rooted Cayley trees. Add
a new vertex $*$ to such a forest and join
$*$ to the root of each tree.
We obtain a Cayley tree with root $*$. This
operation is a one-to-one correspondence, hence 
$Y = q e^Y$. \qed

\begin{corollary}
On the disc $|q| < 1/e$, the function $Y(q)$ is the inverse of
the function $q(Y) = Y/e^Y$.
\end{corollary}

\begin{proposition} \label{Prop:(1-Y)(1+Z)}
We have $(1-Y)(1+Z) = 1$.
\end{proposition}

\paragraph{Proof.} 
$$
Z = DY = D(q e^Y) = q e^Y + q e^Y DY = q e^Y (1+Z)
= Y (1+Z).
$$
Hence $(1-Y)(1+Z) = 1$. \qed

\begin{corollary}
As an abstract algebra, $\A$ is isomorphic to
$\Q [X, X^{-1}]$, where $X=1-Y$.
\end{corollary}

\begin{proposition} \label{Prop:Y^k} We have
$$
Y^k = k \, \sum_{n \geq 1} \frac{(n-1) \dots (n-k+1) n^{n-k}}{n!} \, q^n
= k \sum_{n \geq k} \frac{n^{n-k-1}}{(n-k)!} \, q^n.
$$
\end{proposition}

\paragraph{Proof.} 
Induction on $k$. For $k=1$ the assertion is true. 
To go from $k$ to $k+1$, one uses the equality
$$
D \left(
\frac{Y^{k+1}}{k+1} - \frac{Y^k}k
\right) = (Y^k - Y^{k-1}) DY = Y^{k-1} (Y-1)Z = -Y^k.
$$
It is compatible with our expressions for $Y^k$ and $Y^{k+1}$,
which determines $Y^{k+1}$ up to a constant. But
yhe constant term of $Y^{k+1}$ vanishes. \qed

\bigskip

Now we study the powers of $Z$.

\begin{definition} \label{Def:A_n}
Denote by $A_n$ the sequence of integers
$$
A_n = \sum_{\stackrel{\scriptstyle p+q=n}{p,q \geq 1}}
\frac{n!}{p! \, q!} \, p^p \, q^q.
$$
\end{definition}
Its first terms are $0, 2, 24, 312, 4720, \dots$.
We have
$$
Z^2 = \sum_{n \geq 1} \frac{A_n}{n!} \, q^n.
$$
One can show that 
$$
A_n = n! \, \sum_{k=0}^{n-2} \frac{n^k}{k!}   
\sim \sqrt{\pi/2} \; n^{n+\frac12}.
$$
As far as we know, there is no simple expression for
the powers of $Z$. However, we can prove that they are
linear combinations of the series
$$
D^k Z = \sum_{n \geq 1} \frac{n^{n+k}}{n!} \, q^n
\qquad \mbox{and} \qquad
D^k (Z^2) = \sum_{n \geq 1} \frac{n^k A_n}{n!} \, q^n \, .
$$

\begin{proposition} \label{Prop:Z^k}
For any integer $k \geq 0$, the power series $D^k Z$ and $D^k (Z^2)$
are polynomials in $Z$ with positive integer coefficients,
of degrees  $2k+1$ and $2k+2$ respectively:
\begin{eqnarray*}
D^k(Z) &=& (2k-1)!! \;Z^{2k+1} + \mbox{\rm lower order terms},\\
D^k(Z^2) &=& (2k)!! \; Z^{2k+2} + \mbox{\rm lower order terms}.
\end{eqnarray*}
\end{proposition}

\paragraph{Proof.}
Applying $D$ to both sides of the equality $(1-Y)(1+Z)=1$
we get
$$
-Z(1+Z) + (1-Y) \cdot DZ = 0.
$$
Thus
$$
DZ = \frac{Z(1+Z)}{1-Y} = Z(1+Z)^2.
$$
Hence
$$
D(Z^2) = 2Z^2(1+Z)^2.
$$
Now we proceed by induction on $k$. \qed

\begin{corollary} \label{Cor:Z^k}
For any positive integer $k$, the power series $Z^k$ is a
linear combination with rational coefficients
of the first $k$ series from the list
$Z, Z^2, DZ, D(Z^2), D^2 Z, D^2 (Z^2), \dots$.
\end{corollary}

From Proposition~\ref{Prop:Y^k} and Corollary~\ref{Cor:Z^k}
we deduce the following theorem.

\begin{theorem} {\rm \cite{Zvonkine}}
The algebra $\A$ is spanned over $\Q$ by the power series
$$
1, 
\qquad \quad 
\sum_{n \geq 1} \frac{n^{n+k}}{n!} q^n,
\quad k \in {\mathbb Z}, 
\qquad \quad 
\sum_{n \geq 1} \frac{n^k A_n}{n!} q^n,
\quad k \in {\mathbb N}.
$$
\end{theorem}

Note that the Stirling formula together with the asymptotic
for the sequence $A_n$ allows one to determine the
leading term of the asymptotic for the coefficients of any 
series in $\A$. We have
$$
\frac{n^n}{n!} \sim \frac1{\sqrt{2\pi n}} \; e^n, 
\qquad
\frac{A_n}{n!} \sim \frac12 \, e^n.
$$

Note also that if, for some series $F \in \A$,
we know in advance its degree in $Y$ and in $Z$,
then we can reconstitute the series $F$ using only
a finite number of its initial terms -- a very useful
property for computer experiments.

Combining both remarks, we see that initial terms of the sequence 
of coefficients of $F$ determine the asymptotic of the sequence.

\subsection{Dendrology}
\label{Ssec:dendrology}

\begin{definition} \label{Def:Cayley}
A {\em Cayley tree} is a tree with numbered vertices.
\end{definition}

It is well-known (Cayley theorem) that there are $n^{n-2}$ 
Cayley trees with $n$ vertices. Note that the corresponding
exponential generating function
$$
\sum_{n \geq 1} \frac{n^{n-2}}{n!} \, q^n
$$
lies in the algebra $\A$.

Consider a Cayley tree $T$ with two marked vertices $a$ and $b$. 
Denote by $l(T)$ the distance between these vertices, 
i.e., the number of edges in the shortest path joining them.

\begin{definition}
Denote by $m_{n,k}$ and $p_{n,k}$ the sums 
$$
m_{n,k} = \sum_{T} l(T)^k, \qquad 
p_{n,k} = \sum_{T} \frac{l(T) (l(T)-1) \dots (l(T)-k+1)}{k!}
$$
where the sum is taken over all Cayley trees $T$ with $n$
vertices, two of which are marked.
\end{definition}

For instance, $m_{2,1} = p_{2,1} = 2$. Note that if we consider $l(T)$ as
a random variable, then $m_{n,k}$ is its $k$th moment.

\begin{theorem} \label{Thm:dendrology}
For any $k$, the power series
$$
\sum_{n \geq 1} \frac{m_{n,k}}{n!} \, q^n
\quad \mbox{and} \quad
\sum_{n \geq 1} \frac{p_{n,k}}{n!} \, q^n
$$
lie in $\A$.
\end{theorem}

\begin{example}
It follows from the proof below 
that $p_{n,1} = m_{n,1} = A_n$. This number
is called {\em the total height of Cayley trees} and
was introduced in~\cite{RioSlo}.
\end{example}

\paragraph{Proof of Theorem~\ref{Thm:dendrology}.}
It is sufficient to prove the theorem for $p_{n,k}$.

Fix $k$. There is a natural bijection between the following
sets of objects. 

$E_n$ is the set of Cayley trees with
$n$ vertices, on which one has marked two vertices
by $a$ and $b$ and chosen $k$ distinct edges on the
shortest path from $a$ to $b$. The number of elements
in $E_n$ equals $p_{n,k}$.

$F_n$ is the set of ordered $(k+1)$-tuples of trees with $n$ vertices
in whole; the vertices are numbered from $1$ to $n$ and,
in addition, two vertices $a_i$ and $b_i$, $1 \leq i \leq k+1$,
are marked on each tree.

The bijection is established as follows. Take
a forest from the set $F_n$. Draw new edges $(b_1, a_2)$,
$(b_2, a_3)$, \dots, $(b_k, a_{k+1})$. We obtain a tree with $k$
marked edges lying on the path between $a_1$ and $b_{k+1}$,
i.e., a tree from the set $E_n$.

Now, the trees with two marked vertices are enumerated by
the series $Z$, therefore the exponential generating series
for the sequence $|F_n|$ is $Z^{k+1}$.
\qed

\section{Counting ramified coverings of the sphere}
\label{Sec:coverings}

This section is devoted to the enumeration of 
ramified coverings of the sphere by surfaces of a fixed genus $g$
and to a proof of Theorem~\ref{Thm:Hurwitz}.

\subsection{The ELSV formula}
\label{Ssec:ELSV}

Curiously, the most difficult part of the proof
of Theorem~\ref{Thm:Hurwitz} is the case with only
one multiple ramification point, $k=1$.
We know no other way to prove it
than to use the intersection theory
on moduli spaces. The main ingredient of the
proof is a theorem by 
T.~Ekedahl, S.~K.~Lando, M.~Shapiro, and
A.~Vainshtein that we formulate below after introducing
some notation.

Let $\mu= 1^{a_1} 2^{a_2} \dots$ be a partition
with degeneracy $r$. We define $|\Aut (\mu)|$ to be
$|\Aut (\mu)| = a_1! a_2! \dots$. For the formulation
of the theorem it is more convenient to switch to using
the additive notation for the partition $\mu$,
$\mu = (b_1, \dots, b_p)$, the $b_i$ being the
parts of $\mu$. The Hurwitz number $h_{g,n;\mu}$ is defined in 
Notation~\ref{Not:Hurwitz}.

We denote by $\M_{g,n}$ the {\em moduli space} of smooth genus $g$
curves with $n$ marked and numbered distinct points. 

Further, $\oM_{g,n}$ is the {\em Deligne-Mumford compactification} of
this moduli space; in other words, $\oM_{g,n}$ is the
space of {\em stable} genus $g$ curves with $n$ marked
points.

We denote by $\cL_i$, $1 \leq i \leq n$, the $i$th {\em tautological
line bundle} over $\oM_{g,n}$: consider a point $x \in \oM_{g,n}$
and the corresponding stable curve $C_x$; then the fiber
of $\cL_i$ over $x$ is the cotangent line to $C_x$
at the $i$th marked point. The first Chern class of $\cL_i$
is denoted by $c_1(\cL_i) = \psi_i$.

We will use the expression
$$
\frac{1}{1-\psi_i} = 1 + \psi_i + \psi_i^2 +
\dots \in H^*(\oM_{g,n}, \Q).
$$

Further, we introduce the {\em Hodge vector bundle} $W$
over $\oM_{g,n-r}$. The fiber of
$W$ over a smooth curve is the set of holomorphic $1$-forms
on this curve. The fiber of $W$ over a general stable
curve is the set of global sections of its {\em dualizing sheaf}.
We do not give the details here (see the paper \cite{ELSV}
itself). Suffice it to note that $W$ is a vector bundle
of rank~$g$. 

Now we can write down the ELSV formula.

\begin{theorem} \label{Thm:ELSV} {\rm {\bf (The ELSV formula,}
\cite{ELSV}{\bf )}} \hspace{1em}
For any $g$, $n$, and $\mu$ such that $2-2g-(n-r) <0$, we have
$$
h_{g,n;\mu} \; = \;
\frac{(2n+2g-2-r)!}{|\Aut (\mu)|}
\prod_{i=1}^p \frac{b_i^{b_i}}{b_i!} \; \times
$$
$$
\times \;
\frac1{(n-p-r)!}
\int_{\oM_{g,n-r}} \frac{c(W^*)}
{(1-b_1 \psi_1) \dots (1-b_p\psi_p)
(1-\psi_{p+1}) \dots (1- \psi_{n-r})}
\; .
$$
\end{theorem}

\subsection{Proof of Theorem~\ref{Thm:Hurwitz}}
\label{Ssec:ThmHurwitz}

We are going to prove that all generating series 
$H_{g; \mu_1, \dots,  \mu_k}$ for the Hurwitz numbers
$h_{g,n; \mu_1, \dots, \mu_k}$ (see Notation~\ref{Not:Hurwitz})
with respect to the number
of sheets $n$ lie, once again, in the algebra $\A$.
In Section~\ref{Sec:gravity} we give a motivation for 
considering these particular generating series.

First consider the case of just one partition $k=1$.
This case was essentially covered in~\cite{GoJaVa}. 
Recently, M.~Kazarian~\cite{Kazarian} suggested an improvement of the 
theorem for $k=1$, giving an explicit expression
for the series $H_{g; \mu}$ in terms of the generators $Y$ and $Z$
of $\A$. 

\begin{theorem} {\bf [Kazarian]} \label{Thm:Kazarian}
Consider a partition $\mu = (b_1, \dots, b_p)$, and let $m = \sum b_i$.
Then we have
$$
H_{g;\mu} = \frac{1}{|\Aut(\mu)|} 
\prod_{i=1}^p \frac{b_i^{b_i}}{b_i!} \;
Y^m \; (Z+1)^{2g-2+p} \; \varphi(Z),
$$
where $\varphi(Z)$ is the polynomial
$$
\varphi(Z) = 
\sum_{l \geq 0} \frac{Z^l}{l!}
\int\limits_{\oM_{g,p+l}}
\frac{c(W^*)}{(1-b_1 \psi_1) \dots (1-b_p \psi_p)} \;
\frac{\psi_{p+1}^2 \dots \psi_{p+l}^2}
{(1-\psi_{p+1}) \dots (1-\psi_{p+l})}.
$$
\end{theorem}

Note that the last sum only goes up to $l=3g-3+p$, otherwise
the intergral equals~$0$ for dimension reasons.

\paragraph{Sketch of a proof} (borrowed from~\cite{Kazarian}).
Introduce the series $F = F_{g; b_1, \dots, b_p}$
in an infinite number of variables
$$
F(t_0, t_1, \dots ) = 
\sum_{l; d_1, \dots, d_l}
\frac{t_{d_1} \dots t_{d_l}}{l!}
\int\limits_{\oM_{g,p+l}}
\frac{c(W^*) \; \psi_{p+1}^{d_1} \dots \psi_{p+l}^{d_l}}
{(1-b_1\psi_1) \dots (1-b_p \psi_p)}\; .
$$
We have $H_{g;b_1, \dots, b_p}(q) = q^m F(q,q,q,\dots)$, where
$m = \sum b_i$. On the other hand, $F$ satisfies the string 
and the dilaton equations (see~\cite{Witten}):
\begin{eqnarray*}
\frac{\partial F}{\partial t_0} & = &
m \, F + \sum_{i \geq 0} t_{i+1} \frac{\partial F}{\partial t_i} \; ,\\
\frac{\partial F}{\partial t_1} & = &
\chi \, F + \sum_{i \geq 0} t_i \frac{\partial F}{\partial t_i}\; ,
\qquad \chi = 2g-2+p \; .
\end{eqnarray*}
The theorem now follows from the following fact, obtained by
a manipulation of PDEs.
For any series $F$ satisfying the above string and dilaton equations,
let $\varphi(q) = F(0,0,q,q,q, \dots)$. Then we have
$$
q^m F(q,q,q,\dots) = Y^m \; (1+Z)^\chi \; \varphi(Z).
$$
\qed

\bigskip

Now we deduce the general case from the case $k=1$.

\paragraph{Theorem 2}
{\em Fix any $g \geq 0$, $k \geq 0$. If $g=1$, we suppose that
$k \geq 1$. Then for any partitions
$\mu_1, \dots, \mu_k$, the series
$$
H_{g;\mu_1, \dots, \mu_k} (q) =
\sum_{n \geq 1} \frac{h_{g,n;\mu_1, \dots, \mu_k}}{c(n)!} 
\, q^n
$$
lies in the algebra $\A$.}

\paragraph{Proof of Theorem 2.}
The theorem is proved by induction on the number $k$
of partitions.

{\bf Base of induction.} 
For $k=0,1$, the result is obtained by a direct application
of Theorem~\ref{Thm:Kazarian}. 
In the case $k=0$, we must use Theorem~\ref{Thm:Kazarian}
with an empty partition $\mu$.

There are three exceptional cases in which Theorem~\ref{Thm:Kazarian}
cannot be applied: $g=0$, $k=0$; $g=0$, $k=1$, $p \leq 2$; $g=1$,
$k=0$. These cases are discussed in Remark~\ref{Rem:exceptions}
below. It turns out that the assertion of Theorem~\ref{Thm:Hurwitz}
fails only if $g=1$, $k=0$, as stated in the formulation.

{\bf Step of induction.}
The step of induction is an almost exact repetition of 
the proof of Theorem~2
from our previous work~\cite{Zvonkine}. We only give a short summary of
the argument here.
The proof goes in the spirit of~\cite{GouJac}. 
A similar proof, using the formalizm of colored permutations,
is given in~\cite{Kazarian}.

It is easy to see that there is only a finite number of 
possible cycle structures for a permutation that can be
obtained as a product of two permutations with given cycle
structures $\mu_1$ and $\mu_2$. 

Let $\mu_1$ and $\mu_2$ be two partitions from the list
$\mu_1, \dots, \mu_k$. We can move the two corresponding 
ramification points on $\CP^1$ towards each other until
they collapse. We obtain a new (not necessarily connected)
ramified covering. Its monodromy at the new ramification point
is the product of the monodromies of the two points that have
collapsed. 

Let us choose one
of the possible cycle structures of the product monodromy
and also one of the possible ways in which the 
covering can split into connected components.
By the induction assumption, we obtain a series from
the algebra $\A$ assigned to each connected component of
the covering. Indeed, each connected component is itself
a ramified covering of the sphere as in Theorem~\ref{Thm:Hurwitz},
but with $k-1$ fixed ramification types instead of $k$.
We obtain the generating series for the number of nonconnected 
ramified coverings by multiplying the series that correspond
to the connected components. Since it is a finite product
of series lying in $\A$, we obtain again a series from $\A$.

Finally, we must add the generating series described above
for all choices of types of nonconnected coverings. Since the
number of choices is finite, we obtain, once again,
a series from $\A$.
\qed

\begin{remark} \label{Rem:exceptions}
Let us consider the exceptional cases $g=0$, $k=0,1$
and $g=1$, $k=0$. 

In the genus zero case, the ELSV formula transforms into
a much simpler Hurwitz formula~\cite{Hurwitz,Strehl},
which turns out to be applicable even if the multiple
ramification point has only $1$ or $2$ preimages.

We have, using the notation of Theorem~\ref{Thm:ELSV} and
Notation~\ref{Not:Hurwitz},
$$
h_{0,n;\mu} =
\frac{(2n-2-r)!}{|\Aut(\mu)|} \;
\prod_{i=1}^p \frac{b_i^{b_i}}{b_i!} 
\; \cdot \; 
\frac{n^{n-r-3}}{(n-p-r)!}.
$$
This formula is true for any $n \geq p+r$ and for
any partition $\mu$ (including even the empty partition).
We see that the corresponding generating series always
lies in the algebra $\A$.

The case $g=1$, $k=0$ is covered by the ELSV formula
with an empty partition $\mu$. Consider the moduli 
space $\oM_{1,1}$. Denote by $\beta$ the $2$-cohomology
class of $\oM_{1,1}$ whose integral over the fundamental
homology class equals~1. One can prove that the Hodge
bundle over $\oM_{1,1}$ is a line bundle with first
Chern class $\beta/24$. Therefore we obtain
$$
h_{1,n;\emptyset} = (2n)! \cdot
\frac{1}{n!} \int_{\oM_{1,n}}
\frac{1-\frac{1}{24}\beta}{(1-\psi_1) \dots (1-\psi_n)} \; .
$$
From this we get
$$
\sum_{n \geq 1} \frac{h_{1,n;\emptyset}}{(2n)!} \, q^n
= 
\frac{1}{24} \, \sum_{n \geq 1} \frac{A_n}{n} \frac{q^n}{n!} \; .
$$
This series does not lie in $\A$ (and constitutes the
only exception to the general rule). It suffices
to consider the partition $\mu = (1)$, which
amounts to distinguishing one sheet in the ramified
covering, to obtain the series
$$
\frac{1}{24} \, \sum_{n \geq 1} 
\frac{A_n}{n!} \, q^n \; \in \; \A \, .
$$
\end{remark}

\section{Random metrics and 2-dimensional gravity}
\label{Sec:gravity}

In this section we propose a model of $2$-dimensional gravity
via the enumeration of
ramified coverings. We show that the ``free energy'' function
and the values of ``observables'' coinside with those
obtained in other models.

We also draw a parallel between the
study of spaces of Riemannian metrics using
ramified coverings of the sphere and the study of
spaces of abelian differentials using ramified coverings
of the torus.

\subsection{Models of 2-dimensional gravity}
\label{SSec:2models}

Here we explain what sort of questions about ramified coverings
arise in $2$-dimensional gravity and why. Precise mathematical
results are given below in Sections \ref{Ssec:Painleve} 
and~\ref{Ssec:KdV}. 

In every problem of statistical physics one starts
with introducing a space of states and by assigning
an energy to every state. 

In 2-dimensional gravity, a {\em state} is a
2-dimensional compact oriented real
not necessarily connected surface endowed with a
Riemannian metric. Two surfaces like that are
equivalent, i.e., correspond to the same state,
if they are isometric.

Consider a surface $S$ with a Riemannian metric. Let
$\chi(S)$ be its Euler characteristic and $A$ its
total area. To such a surface one assigns an {\em energy}
$$
E = \lambda A + \mu \chi(S).
$$
Here $\lambda$ and $\mu$ are two constants called
the cosmological constant and the gravitational
constant, respectively. Note that $\chi(S)$ is
actually the integral over $S$ of the scalar curvature
of the metric. The fact that this integral takes
such a simple form is special to dimension~2.

Now the first thing to do is to compute the
{\em partition function} $z(\lambda, \mu)$
or, equivalently, the {\em free energy} $f(\lambda,\mu)$
$$
z(\lambda, \mu) = \int_{\mbox{states}} e^{-E},
\qquad \qquad
f(\lambda, \mu) = \ln z(\lambda, \mu)
= \sum_{g \geq 0} \int_{\mbox{metrics}}
e^{-E}.
$$
The free energy is the sum of contributions of
connected surfaces, while the partition function
is the sum of contributions of all surfaces.

Neither of the above integrals is well-defined
mathematically, but we would still like to
compute them. To do that, physicists introduced
a {\em discrete model} of Riemannian metrics,
replacing them by quadrangulations~\cite{BreKaz,Witten}
(see also~\cite{thebook}, Chapter~3 
for a mathematical description). In this model,
instead of considering Riemannian metrics, one considers
metrics obtained by gluings of squares of area $\varepsilon$.
Our goal is to show that the ramified coverings
of the sphere provide a new (maybe more natural)
discrete model of Riemannian metrics.

Fix a positive number $\varepsilon$. Consider a sphere
with the standard (round) Riemannian metric
of total area $\varepsilon$. On this sphere,
choose at random $2n+2g-2$ points. Now chose
a random connected $n$-sheeted covering of the sphere with simple
ramifications over the $2n+2g-2$ chosen points. 
The covering surface $S$ will automatically be of genus $g$. 
The metric on the sphere can be lifted to $S$, which will
give us a metric with constant positive curvature except
at the critical points, where it has conical singularities
with angles $4\pi$. This metric is, of course, not Riemannian.
However, one can argue that if $\varepsilon$ is very small and the number
of sheets very large, a random metric obtained
in this way looks similar to a random
Riemannian metric (unless we look at them through a
microscope to reveal the difference). We do not know
any rigorous statement that would formalize this
intuitive explanation, but the same argument is
used by physicists to justify the usage of quadrangulations.

Using our discrete model of metrics, one can
write the free energy for the 2-dimensional
gravity in the following way:
$$
f(\lambda, \mu) = 
\sum_{g,n} \frac{\varepsilon^{2n+2g-2}}{(2n+2g-2)!} \; 
h_{g,n; \emptyset} \; e^{-\lambda n \varepsilon - \mu (2-2g)}.
$$
Here the factor $\varepsilon^{2n+2g-2}/(2n+2g-2)!$ is the volume of the space
of choices of $2n+2g-2$ unordered points on the sphere
of area $\varepsilon$, while $n \varepsilon$ in the exponent
is the area of the covering surface.

In Section~\ref{Ssec:Painleve} we show that
if we let $n \rightarrow \infty$ while $g$ remains
fixed, we have
$$
\frac{h_{g,n;\emptyset}}{(2n+2g-2)!} \sim e^n \; n^{\frac52 (g-1)-1}
\; b_g
$$
for some constants $b_g$. Thus the coefficients of $f$ have the 
following asymptotic:
$$
\frac{\varepsilon^{2n+2g-2}}{(2n+2g-2)!} \; 
h_{g,n; \emptyset} \; e^{-\lambda n \varepsilon - \mu (2-2g)}
\; \sim \; b_g \; 
e^{- (\lambda \varepsilon - 2 \ln \varepsilon +1)n} \;
(\varepsilon e^\mu)^{2g-2} \; n^{\frac52(g-1) -1}.
$$

Now we make the final step by letting $\varepsilon$ tend
to $0$ in the expression of~$f$. 
To obtain an interesting limit for the free energy, 
we must {\em make $\lambda$ and $\mu$ depend on $\varepsilon$}.
We want to use
$$
\sum_{n \geq 1} n^{\gamma-1} e^{-\delta n}  \sim 
\frac{\Gamma(\gamma)}{\delta^\gamma}
\quad \mbox{ as } \delta \rightarrow 0.
$$
Therefore we set $\gamma = \frac52 (g-1)$ and we let
$\delta = \lambda \varepsilon -2 \ln \varepsilon - 1$ tend to~$0$, 
while
$$
y = 
\frac{(\varepsilon e^\mu)^{4/5}}\delta
=
\frac{(\varepsilon e^\mu)^{4/5}}
{\lambda \varepsilon -2 \ln \varepsilon-1}
$$
remains fixed. This gives us the final expression
of the free energy, now depending on only one variable
$y$:
$$
f(y) = \Gamma (-5/2) \, b_0 \, y^{5/2} - b_1 \ln y
+ \sum_{g \geq 2} 
\Gamma \biggl(5(g-1)/2\biggr) \,
\, b_g \cdot y^{5(1-g)/2}.
$$
The coefficients
$\Gamma (5(g-1)/2) \, b_g $
are rational for odd $g$ and rational multiples
of $\sqrt2$ for even $g$.

\bigskip

Our above treatment is parallel to E.~Witten's treatment of
the quadrangulation model in~\cite{Witten}. 
Denote by $Q_{g,n}$ the number of ways to
divide a surface of genus $g$ into $n$ squares.
Then the study of the quadrangulation model
involves the asymptotic of $Q_{g,n}$, which is
given by
$$
Q_{g,n} \sim 12^n n^{\frac52(g-1) -1} b'_g,
$$
for another sequence of constants $b'_g$. This
sequence was studied using matrix integrals, and it
is known that a generating function for the sequence
$b'_g$ satisfies the Painlev{\'e} I equation. In the next
section we show a similar result for the constants $b_g$.
This implies that the functions $f$ obtained in the two
models coincide up to a rescaling of the variable $y$;
more precisely, we have
$$
b_g' = 2^{\frac32(g-1)+1} b_g \, .
$$

\bigskip

In the treatment of the quadrangulation model in~\cite{Witten},
Witten also introduced {\em observables} that correspond to
counting quadrangulations with ``impurities'', that is,
the number of ways to divide a surface into a large number
of squares and a fixed number of given polygons. Each 
observable $\tau_d$ is represented by a formal linear combination
of a $2$-gon, a $4$-gon, and so on, up to a $(2d+2)$-gon.
The values of these observables combine into a generating
function $F(t_0, t_1, \dots)$ that can be studied using matrix integrals.
It turns out that $\partial^2 F /\partial t_0^2$ is
a solution of the Korteweg--de~Vries (KdV)
hierarchy. This solution is called the ``string solution''.

From now on the notation $\left< \cdot \right>$ 
or $\left< \cdot \right>_g$ will mean
$$
\left< \tau_{d_1} \dots \tau_{d_n} \right>_g = 
\left< \tau_{d_1} \dots \tau_{d_n} \right> 
= \int\limits_{\oM_{g,n}}
\psi_1^{d_1} \dots \psi_n^{d_n}.
$$
It turns out that the generating series
$$
F(t_0, t_1, \dots)
= \sum_{n \geq 1} \sum_{d_1,\dots,d_n}
\left< \tau_{d_1} \dots \tau_{d_n} \right> 
\frac{t_{d_1} \dots t_{d_n}}{n!} \; .
$$
coincides with the series $F$ obtained from the quadrangulation
model. In particular, its second derivative 
$U = \partial^2 F / \partial t_0^2$. 
satisfies the KdV equation:
\begin{equation}\label{Eq:KdV}
\frac{\partial U}{\partial t_1} =
U \frac{\partial U}{\partial t_0} +
\frac1{12} \frac{\partial^3 U}{\partial t_0^3}.
\end{equation}
This was conjectured by E.~Witten in~\cite{Witten} and 
proved by M.~Kon\-tsevich in~\cite{Kontsevich}. 

In Section~\ref{Ssec:KdV} we show that the numbers
$\left< \tau_{d_1} \dots \tau_{d_n} \right>$ can
also be obtained in the model of ramified coverings.
Each observable $\tau_d$ is represented by a formal
linear combination of a noncritical point, a simple
critical point, and so on, up to a $d$-tuple critical
point. 

\bigskip

Recently, M.~Kazarian and S.~K.~Lando~\cite{KazLan}
found an independent proof of the fact that the function
$U$ arising in the enumeration of Hurwitz numbers 
satisfies the KdV equation. This has lead them to
a new proof of Witten's conjecture.

\subsection{The Painlev{\'e} I equation}
\label{Ssec:Painleve}

The results of this section were obtained in common with S.~Lando.

\begin{proposition}
For a fixed $g$, we have
$$
\frac{h_{g,n; \emptyset}}{(2n+2g-2)!} 
\sim e^n n^{\frac52 (g-1) - 1} b_g \qquad \mbox{as} 
\quad n \rightarrow \infty,
$$
where
$$
b_g = \frac{ \left< \tau_2^{3g-3} \right> }
{(3g-3)! \; \; 2^{\frac52(g-1)} \; \Gamma \left(\frac52(g-1) \right)}
\; .
$$
\end{proposition}

\paragraph{Proof.} From Theorem~\ref{Thm:Kazarian} we see
that $H_{g;\emptyset}$ is a polynomial in $Z$ with
leading term
$$
\frac{ \left< \tau_2^{3g-3} \right> } {(3g-3)!} \; Z^{5g-5}. 
$$
Indeed, if $l=3g-3$, the degree of the class $\psi_1^2 \dots \psi_l^2$
is exactly the dimension of $\oM_{g,l}$ (both are equal to $6g-6$).
Therefore the classes $c(W^*)$ and $1/(1-\psi_i)$ do
not contribute. On the other hand, by Proposition~\ref{Prop:Z^k}
the coefficients of the series $Z^l$ grow as
$$
\frac{n^{\frac{l}2-1} \; e^n}{\Gamma(l/2) \; 2^{l/2}}.
$$
Multiplying the coefficient of the leading term
$$
\frac{\left< \tau_2^{3g-3} \right>}{(3g-3)!} Z^{5g-5}
$$ 
of $H_g$ by the asymptotic of coefficients of $Z^{5g-5}$
we obtain the leading term of the asymptotic of $h_{g,n}$,
in particular, the constant $b_g$. \qed

The first values of the constants $b_g$ are
$$
b_0 =\frac1{\sqrt{2\pi}}, \quad
b_1 =\frac{1}{2^4 \cdot 3}, \quad
b_2 = \frac1{\sqrt{2\pi}} \; \frac{7}{2^5 \cdot 3^3 \cdot 5},
$$
$$
b_3 = \frac{5 \cdot 7^2}{2^{16} \cdot 3^5}, \quad
b_4= \frac1{\sqrt{2\pi}} \;
\frac{7 \cdot 5297}{2^{11} \cdot 3^8 \cdot 5^2 \cdot 11 \cdot 13}.
$$
The above expression for $b_g$ allows us to rewrite the
function $f''(y)=u(y)$ in the following way:
\begin{equation} \label{Eq:2}
u(y) = -\sqrt{2y} + \frac1{12} (2y)^{-2}
+ \sum_{g \geq 2} (5-5g)(3-5g) 
\frac{\left< \tau_2^{3g-3} \right>}{(3g-3)!} \; 
(2y)^{\frac{1-5g}{2}} \; .
\end{equation}

\begin{proposition} \label{Prop:Painleve}
The function $u(y)$ satisfies the
Painlev{\'e}~I equation
$$
u^2(y) + \frac16 u''(y) = 2y.
$$
\end{proposition}

\paragraph{Proof.} 

By extracting the coefficient of $t_2^{3g-1}$ in
the KdV equation~(\ref{Eq:KdV}), we obtain, for every $g \geq 1$,
\begin{equation} \label{Eq:3}
\frac{\left< \tau_0^2 \tau_1 \tau_2^{3g-1} \right>_g}{(3g-1)!} \;
= \!\!\!\!\!
\sum_{\stackrel{\scriptstyle g'+g''=g}{g' \geq 1, g'' \geq 0}}
\!\!\!\!
\frac{\left< \tau_0^2 \tau_2^{3g'-1} \right>_{g'}}{(3g'-1)!}
\cdot
\frac{\left< \tau_0^3 \tau_2^{3g''} \right>_{g''}}{(3g'')!} 
\; + \;
\frac1{12} \,
\frac{\left< \tau_0^5 \tau_2^{3g-1} \right>_{g-1}}{(3g-1)!}  \, .
\end{equation}

Now we use the string and the dilaton equations
to kill the $\tau_0$ and $\tau_1$ factors in all the brackets
of~(\ref{Eq:3}). We obtain the following identities (with some
exceptions for low genus):

\begin{eqnarray*}
\frac{\left< \tau_0^2 \tau_1 \tau_2^{3g-1} \right>_g}{(3g-1)!}
&=&
 (5g-5)(5g-3)(5g-1) 
\frac{\left< \tau_2^{3g-3} \right>_g}{(3g-3)!} \; , 
\qquad \;\; \frac{\left< \tau_0^2 \tau_1 \tau_2^2 \right>_1}{2!} = 
\frac13 \; . \\
\\
\frac{\left< \tau_0^2 \tau_2^{3g'-1} \right>_{g'}}{(3g'-1)!}
&=& 
(5g-5)(5g-3) 
\frac{\left< \tau_2^{3g'-3} \right>_{g'}}{(3g'-3)!} \; , 
\qquad \qquad \qquad \; 
\frac{\left< \tau_0^2 \tau_2^2 \right>_1}{2!} = 
\frac1{12} \; . \\
\\
\frac{\left< \tau_0^3 \tau_2^{3g''} \right>_{g''}}{(3g'')!}
&=& 
(5g-5)(5g-3)(5g-1) 
\frac{\left< \tau_2^{3g''-3} \right>_{g''}}{(3g''-3)!} \; , \\
& & \left< \tau_0^3 \right>_0= 1 \, ,
\quad 
\frac{\left< \tau_0^3 \tau_2^3 \right>_1}{3!} 
= \frac13 \; . \\
\\
\frac{\left< \tau_0^5 \tau_2^{3g-1} \right>_{g-1}}{(3g-1)!}
&=& 
(5g-10)(5g-8)(5g-6)(5g-4)(5g-2) 
\frac{\left< \tau_2^{3g-6} \right>_{g-1}}{(3g-6)!} \; ,\\
& & \frac{\left< \tau_0^5 \tau_2^2 \right>_0}{2!}= 3 \, ,
\quad 
\frac{\left< \tau_0^5 \tau_2^5 \right>_1}{5!} = 16 \, . \\
\end{eqnarray*}
Substitute these expressions in~(\ref{Eq:3}) and
compare to the expression~(\ref{Eq:2}) of $u$.
Taking into account the exceptional starting terms
we obtain
$$
\sqrt{2y} \left( u'(y)+\frac1{\sqrt{2y}} \right)
= \left( u(y)+\sqrt{2y} \right) \; u'(y) + \frac1{12} u'''(y).
$$
We rewrite this as
$$
u(y) \, u'(y) + \frac1{12} u'''(y) = 1
$$
and integrate it once to obtain
$$
u^2(y) + \frac1{6} u''(y) = 2y.
$$
\qed

\subsection{The KdV hierarchy}
\label{Ssec:KdV}

\setlength\unitlength{1mm}

Now we show how to obtain the numbers
$\left< \tau_{d_1} \dots \tau_{d_p} \right>$
as leading term coefficients of the asymptotic of Hurwitz numbers.

Note that Okounkov and Pandharipande~\cite{OkoPan} also
obtained the numbers $\left< \tau_{d_1} \dots \tau_{d_p} \right>$
using asymptotics of Hurwitz numbers. However their asymptotics
are different from ours and have no direct physical interpretation.

The Hurwitz numbers involved are those that count
ramified coverings with many simple
ramification points, but {\em only one}
multiple ramification point that has
$p$ marked preimages. Each marked preimage 
corresponds to a factor $\tau$ in the bracket. Using
the notation from Section~\ref{Sec:coverings}, the numbers
we are interseted in are $|\Aut (\mu)| \cdot  h_{g,n; \mu}$, with
$\mu = (b_1, \dots, b_p)$. The factor $|\Aut (\mu) |$ is due
to the fact that the marked preimages are numbered.

Denote by $\fbox{b}$ a $b$-tuple preimage of the special
ramification point.
Our result is then best described by the following symblic formula:
$$
\tau_d  \; \mbox{``}=\mbox{''} \; 
\frac1{0! \, (d+1)^d} \cdot
\frac{\fbox{$d+1$}}{Y^{d+1}}
-
\frac1{1! \, d^{d-1}} \cdot
\frac{\fbox{$d$}}{Y^{d}}
+ \dots
+ (-1)^d
\frac1{d! \, 1^0} \cdot
\frac{\fbox{$1$}}{Y}
$$
The recipe for obtaining the number 
$\left< \tau_{d_1} \dots \tau_{d_p} \right>$
is the following.

\noindent
{\bf 1.}~~Replace each $\tau_d$ by the right-hand side of the above
symblic equality. 

\noindent
{\bf 2.}~~Expand the product to obtain a linear combination of
terms of the form
$$
\mbox{const} \cdot \frac{\fbox{$b_1$} \dots \fbox{$b_p$}}
{Y^{b_1+ \dots + b_p}} \; .
$$

\noindent
{\bf 3.}~~To each such term assign the partition
$\mu = (b_1, \dots, b_p)$ and the series
$$
\mbox{const} \cdot \frac{|\Aut(\mu)| \, H_{g;\mu}}{Y^{b_1 + \dots + b_p}}.
$$

\noindent
{\bf 4.}~~Add all the series thus obtained. This gives a series in
$\A$ that we denote by
$H[\tau_{d_1} \dots \tau_{d_p}] (q).$

\begin{theorem}
We have
$$
H[\tau_{d_1} \dots \tau_{d_p}] \; = \;
\left< \tau_{d_1} \dots \tau_{d_p} \right>
\; (Z+1)^{2g-2+l}.
$$
The asymptotic of the coefficient of $q^n$ 
{\rm (}as $n \rightarrow \infty${\rm )}
in $H[\tau_{d_1} \dots \tau_{d_p}] (q)$ equals
$$
\frac{\left< \tau_{d_1} \dots \tau_{d_p} \right>}
{2^{\frac{2g-2+p}{2}} \; \Gamma\left(\frac{2g-2+p}{2}\right)}
\cdot 
e^n \; n^{\frac{2g-2+p}{2}-1} \; .
$$
\end{theorem}

\paragraph{Proof.} 
Here again we will use Theorem~\ref{Thm:Kazarian}.
The crucial part of this proposition is the polynomial 
$$
\varphi_{b_1, \dots, b_p} (Z) = 
\sum_{l \geq 0} \frac{Z^l}{l!}
\int\limits_{\oM_{g,p+l}}
\frac{1}{(1-b_1 \psi_1) \dots (1-b_p \psi_p)} \;
\frac{c(W^*) \; \psi_{p+1}^2 \dots \psi_{p+l}^2}
{(1-\psi_{p+1}) \dots (1-\psi_{p+l})}.
$$
We are going to consider linear combinations
of such polynomials for different $b_i$'s. Our goal
is to obtain a cancellation of all higher order terms
in $Z$ leaving only a constant term ($l=0$). This constant term
will turn out to be 
$\left< \tau_{d_1} \dots \tau_{d_p} \right>$.

First of all, here is a linear combination of the
series $1/(1 - \psi)$, $1/(1-2\psi)$, \dots, $1/(1-(d+1) \psi)$
whose terms up to $\psi^{d-1}$ vanish:
\begin{equation} \label{Eq:lincomb}
\frac{1}{d!}
\sum_{b=1}^{d+1} 
\frac{(-1)^{d+1-b} \, {d \choose b-1}}{1-b \psi}
= \psi^d + O(\psi^{d+1}).
\end{equation}

In the expression
$$
|\Aut(\mu)| \, H_{g;\mu}
= (Z+1)^{2g-2+p} \; Y^{b_1 + \dots + b_p}
\prod_{i=1}^p \frac{b_i^{b_i}}{b_i!} \;  \cdot \; 
\varphi(Z),
$$
the integrand $1/(1-b_i \psi_i)$ appears with an additional
factor $Y^{b_i} b_i^{b_i}/b_i!$ in front of the
integral. To compensate for this factor, we multiply
the coefficients of~(\ref{Eq:lincomb}) by its inverse,
which gives
$$
\sum_{b=1}^{d+1} 
\frac{(-1)^{d+1-b}}{(d+1-b)!\; b^{b-1}}
\cdot \frac{\fbox{b}}{Y^b}.
$$
This is precisely the formula that we gave for $\tau_d$.

Multiplying such expressions for $d=d_1, \dots, d_p$
and adding them up we obtain
$$
H[\tau_{d_1} \dots \tau_{d_n}] =
$$
$$
(Z+1)^{2g-2+p} \cdot
\sum_{l \geq 0} \frac{Z^l}{l!}
\int\limits_{\oM_{g,p+l}}
\bigl(
\psi_1^{d_1} \dots \psi_p^{d_p} + \mbox{h.o.t.}
\bigr) \; 
\frac{c(W^*)\; \psi_{p+1}^2 \dots \psi_{p+l}^2}
{(1-\psi_{p+1}) \dots (1-\psi_{p+l}) }
\;,
$$
where ``h.o.t.'' means ``higher order terms''.

The above integral vanishes for dimension reasons whenever $l>0$.
For $l=0$, the factor $c(W^*)$ contributes only by $c_0(W^*) = 1$.
Thus $H[\tau_{d_1} \dots \tau_{d_n}] = 
\left< \tau_{d_1} \dots \tau_{d_p} \right> \cdot (Z+1)^{2g-2+p}$
as claimed.

The second assertion of the theorem follows from the first one
and from the asymptotic of the coefficients of $Z^k$
(Proposition~\ref{Prop:Z^k}). \qed

\begin{remark} [Kazarian]
Equality~(\ref{Eq:lincomb}) can be used to express the numbers
$\left< \tau_{d_1} \dots \tau_{d_n} \right>$ as finite
linear combinations of Hurwitz numbers, without considering
any asymptotics. More precisely, we have
$$
\left< \tau_{d_1} \dots \tau_{d_n} \right> =
\sum_{b_1, \dots, b_p} \; 
\left( \prod_{i=1}^p \frac{(-1)^{d_i+1-b_i}}{(d_i+1-b_i)! \,
    b_i^{b_i-1}}
\right)
\frac{|\Aut(b_1, \dots, b_p)| \; h_{g,n; b_1, \dots, b_p}}
{c(n)!}.
$$
Here the sum is over $1 \leq b_i \leq d_i+1$, the number
of sheets is $n = \sum b_i$ and $c(n) = n+p+2g-2$.
\end{remark}

\subsection{Ramified coverings of a torus and abelian
differentials}

Fix an integer $g \geq 1$ and a list of $p$ nonnegative
integers $b_1, \dots, b_p$ with the condition $\sum b_i = 2g-2$. 
We consider the space
$D_{g;b_1, \dots, b_p}$ of abelian differentials
on Riemann surfaces of genus $g$, with zeroes of
multiplicities $b_1, \dots, b_p$. More precisely,
$D_{g;b_1, \dots, b_p}$ is the space of triples
$(C, \{ x_1, \dots, x_p \} , \alpha)$, where
$C$ is a smooth complex curve, $x_1, \dots, x_p \in C$
are distinct marked points, and $\alpha$ is an
abelian ($=$ holomorphic) differential on $C$ whose
zero divisor is precisely $b_1 x_1 + \dots + b_p x_p$.

It turns out that the space $D_{g;b_1, \dots, b_p}$
has a natural {\em integer affine structure}.
This means that it can be covered by charts of
local coordinates in such a way that the transition functions
are affine maps with integer coefficients. Such 
local coordinates are introduced as follows. 
Fix a basis $l_1, \dots, l_{2g+p-1}$ of 
the relative homology group
$H_1(C,\{x_1, \dots, x_p \}, {\mathbb Z})$.
Then the integrals of $\alpha$ over the cycles
$l_i$ are the local coordinates we need.
The {\em area function}
$$
A: (C, \alpha) \mapsto \frac{i}{2}
\int_C \alpha \wedge {\bar \alpha}
$$
is a quadratic form with respect to the affine
structure.

The integer affine structure allows one to define
a volume measure on the space
$D_{g;b_1, \dots, b_p}$. It is then a natural
question to find the total volume of the part
of the space $D_{g;b_1, \dots, b_p}$ defined
by $A \leq 1$ (the volume of the whole space
being infinite).

A.~Eskin and A.~Okounkov~\cite{EskOko} obtained an
effective way to calculate these volumes using
the asymptotic for the number of ramified coverings
of a torus. Consider the elliptic curve obtained
by gluing the opposite sides of the square
$(0,1,i,1+i)$ endowed with the abelian differential
$dz$. Given a ramified covering of this elliptic curve
with critical points of multiplicities
$b_1, \dots, b_p$, we can lift the abelian differential
to the covering curve and obtain a point of
$D_{g;b_1, \dots, b_p}$. One can then easily show that 
such points are densely and uniformly distributed in
$D_{g;b_1, \dots, b_p}$ if one considers coverings
with a big number of sheets. Moreover, 
R.~Dijkgraaf~\cite{Dijkgraaf} and S.~Bloch and
A.~Okounkov~\cite{BloOko} showed that the generating series for
the ramified coverings of the torus that arise
in this study are quasi-modular forms. In other
words, they lie in the algebra
$$
\Q [ E_2, E_4, E_6],
$$
where $E_{2k}$ are the Eisenstein series
$$
E_{2k}(q) = \frac12 \zeta(1-2k) + \sum_{n \geq 1}
\left( \sum_{d | n} d^{2k-1} \right) q^n.
$$

We conclude with the following comparison between the
counting of ramified coverings of a sphere and of a torus.

\bigskip 

Sphere: The generating series enumerating the ramified
coverings lie in the algebra $\A$.

Torus: The generating series enumerating the ramified
coverings lie in the algebra of quasi-modular forms.

\bigskip

Sphere: The coefficients of a
generating series grow as $e^n \cdot n^{\gamma-1} \cdot c$.
The exponent $\gamma$ is a half-integer. The
number $- \gamma$ is called the
{\em string susceptiblity}.
The constant $c$ is an observable in 2-dimensional gravity.

Torus: The  sum of the first $n$ coefficients of a
generating series grows as $n^d \cdot c$. The number
$d$ is the complex dimension of the corresponding space of
abelian differentials. The constant $c$ is its volume.

\bigskip

Sphere: The observables can be arranged into a generating
series that is a solution of the KdV hierarchy.

Torus: As far as we know, nobody has tried to arrange the
volumes of the spaces of abelian differentials into
a unique generating series.

\end{document}